\documentclass[letterpaper, 10 pt, conference]{ieeeconf}  

\usepackage{cite}
\usepackage{caption}
\usepackage{subcaption}
\usepackage{hyperref}
\usepackage{amsmath}
\usepackage{amsmath,amssymb,amsfonts}
\usepackage{algorithmic}
\usepackage{graphicx}
\usepackage{float}
\usepackage{textcomp}
\usepackage{xcolor}
\title{\LARGE \bf
Drone Delivery Optimization
}

\author{Saayuj Deshpande$^{1}$ Purushotham Mani$^{1}$  
}

\begin{document}

\maketitle
\thispagestyle{empty}
\pagestyle{empty}

\begin{abstract}

This research has addressed three critical challenges inherent in the implementation of drone delivery systems, namely, optimizing battery charging station placement, solving the shortest path problem for drones within their single battery charge travel distance, and efficiently scheduling multiple drones across numerous warehouses and delivery locations with diverse demands. The study has leveraged a 2D grid model with obstacles, providing a practical foundation extendable to a 3D grid for accommodating complex structures. For battery station placement, the Miller-Tucker-Zemlin subtour elimination method has been applied to avoid the formation of charging station clusters. Future research directions involve the integration of these cases into a holistic solution, exploration of three-dimensional space, and the pursuit of bi-level optimization considering the interdependence of battery station placement and shortest path determination. This study contributes to the emerging field of drone delivery systems by addressing key optimization challenges and paving the way for comprehensive, integrated solutions.

\textit{Keywords}--- Drone delivery systems, Battery charging station, Shortest path problem, Multi-drone scheduling, Mixed integer programming, Miller-Tucker-Zemlin subtour elimination

\end{abstract}

\section{Introduction}

E-commerce has led to increasing demand for efficient and cost-effective last-mile delivery services. Drones in delivery services have gained significant attention in recent years due to their potential to reduce delivery times and costs. Furthermore, drones are less expensive to maintain than traditional delivery vehicles and have the ability to traverse terrains inaccessible by road. Although helicopters are widely used in emergency situations to transport food and medical supplies, drones are a much cheaper alternative. The use of unmanned vehicles like drones also reduces the labour cost significantly for private industries. They do have their drawbacks though. Their flight time and payload capacity are limited, which puts them at a major disadvantage against conventional delivery vehicles. Maintaining drone systems demands expertise surpassing that needed for traditional transportation. The complexities of drone technology, encompassing unmanned aerial vehicle (UAV) systems, real-time analytics, and regulatory navigation, require a highly skilled workforce.

The efficient deployment of drones requires addressing several challenges, such as the drone’s single-charge travel distance, placement of battery charging stations, and ensuring collision avoidance among drones and obstacles.

This project has addressed these challenges by developing a drone delivery system that optimizes the location of the battery charging stations, generates the shortest delivery route in the modeled grid, and handles the scheduling of multiple drones operating across multiple delivery locations and warehouses with varied demand to reduce overall delivery cost and time.

The geographical area has been modeled as a two-dimensional grid, with blocked cells or obstacles representing restricted areas where drone travel is prohibited. In contrast, unblocked cells represent accessible regions where the drone can travel freely, and where battery charging stations can be placed. The model has been formulated using both linear and non-linear constraints and has been solved using the Gurobi solver within the AMPL programming framework and Python. It is believed that this approach has significantly improved the efficiency of drone delivery services, and can be further extended to address various other challenges in this domain.

\section{Related Work}

While extensive research specific to the problems addressed in this paper has been limited, it is essential to provide an overview of the existing work within this domain. \cite{8891117} has formulated a model with multiple objectives: to minimize the total delivery cost, to minimize the total number of unsuccessful delivered packages, and to maximize the reward of on-time delivery. Furthermore, this optimization has been structured as a three-stage stochastic programming approach, designed to address operational uncertainties. The drone delivery problem has been approached from a more commercial angle in \cite{7513397}. This research has addressed the optimization of drone delivery with a focus on minimizing the delivery time and operational costs. A key factor considered in this study has been the power consumption of the multi-rotor copters during their operations. The research has leveraged both Mixed Integer Linear Programming and a simulated annealing algorithm to derive its results. \cite{7934790} has proposed an operation management strategy based on the optimization of modular drone delivery systems. The module has been defined on the basis of the functionality of its components that define its type and a set of variables that differentiates variants of a module type.

Aspects of the Traveling Salesman Problem (TSP) have been utilized to address the drone delivery problem in the studies cited as \cite{article1}, \cite{Ha_2019} and \cite{HA2018597}. These variants are popularly referred to as TSP-D problems. Research has also been done to solve the TSP-D problem using the popular dynamic programming concept, as in \cite{articletspddp}. Some aspects of TSP-D have been utilized in our first case, which involves the optimal placement of battery charging stations. An essential step to carry out during any TSP or Vehicle Routing Problem (VRP) is subtour elimination. \cite{doi:10.1057/jors.1986.86} has discussed several subtour elimination methods, most of which have employed cardinality constraints on the subset of nodes, similar to the Dantzig-Fulkerson-Jhonson method. The Miller-Tucker-Zemlin method has been discussed in detail in \cite{DESROCHERS199127}, and improvements have been suggested to strengthen these constraints and generalize them to a variety of VRPs. \cite{ACHUTHAN1996573} has described a new subtour elimination constraint exploiting the cardinality constraints themselves, while \cite{Pferschy_2016} has formulated constraints specifically for pure integer solutions, which considerably reduced the computation time.

To aid with the distance calculation between different nodes on our grid, Dijkstra's algorithm has been utilized. In \cite{articledijk}, the application of Dijkstra's algorithm for finding the shortest path between buildings has been demonstrated. The algorithm simply needs a graph, including its vertices, edges, and weights, to be given as input. \cite{5359145} has described a new shortest path algorithm based on Dijkstra, which has been shown to work faster on large datasets. Applications of Dijkstra to the route planning problem have been discussed in \cite{5569452} and \cite{ApplyingDijkstrasAlgorithminRoutingProcess}.

A pivotal aspect of our research has focused on addressing the significant challenge posed by battery limitations in drone operations. In \cite{article2}, a route planning model has been presented, which has factored in the battery consumption rate induced by the payload. Alternatively, another approach to account for battery limitations has involved establishing a clear relationship between payload and drone flight range, defined in terms of battery capacity, as exemplified in \cite{inproceedings}. \cite{article3}, \cite{article4}, and \cite{9068270} have dealt with the optimal deployment of battery charging stations within a specified demand area. These studies have taken into account a multitude of factors, including drone flight range, the coverage of recharging stations for delivery services, and the establishment of a viable delivery network comprising warehouses and battery charging stations. These strategies have offered valuable insights into how the critical issue of managing battery constraints has been tackled within our research context.

An issue that has not been tackled in our paper pertains to the consideration of challenges arising from drone reliability. In reference \cite{8453380}, a model has been outlined that provides an estimation of the expected demand loss attributable to drone failures. In \cite{9655791}, the underlying algorithm used to solve the problem is the same as the Traveling Salesman Problem but it has also delved into the commercial dimensions by incorporating delivery authentication mechanisms, such as QR code verification, facilitated through image processing libraries like OpenCV. Both of these subjects have exhibited strong applicability to real-world scenarios, emphasizing the necessity for further investigation.

\section{Problem Setting}

\subsection{Optimizing the Location of Battery Charging Stations}
This research has tackled an optimization problem centered on the efficient placement of battery charging stations within a grid framework featuring obstacles. The core objective is to minimize the number of charging stations while ensuring comprehensive grid coverage, crucial for enabling the effective operation of drones with restricted travel distances per battery charge.

To address this problem, the requirement emerged for determining the shortest distance between each node and every other node within our grid. However, the presence of obstacles rendered a straightforward Euclidean distance calculation inadequate. To overcome this challenge, Dijkstra's algorithm has been utilized for distance computation, assigning edges between each node and its neighbors a distance of infinity in the event that the neighbor represents an obstacle, and the corresponding Euclidean distance otherwise. This has enabled the computation of the travel distance between any two points on the grid.

Once the solution demonstrated the placement of battery charging stations ensuring the reachability of every node, the subsequent challenge was guaranteeing interconnectedness among the stations. Employing a basic flow conservation approach led to the formation of clusters, each internally reachable but isolated from stations in other clusters. In response, the Miller-Tucker-Zemlin (MTZ) subtour elimination technique has been implemented to address and resolve this issue. In our MTZ implementation, a depot battery charging station has been randomly placed in the interior and assigned a rank of 1. Subsequently, each additional battery charging station has received a monotonically increasing rank, accompanied by a constraint stipulating that any outgoing edge from a battery charging station must connect to a higher-ranked station. This constraint has effectively prevented the formation of subtours, contributing to the full realization of our study's objective.

\subsection{Shortest Path Problem}
This study has investigated the problem of path planning for a drone operating within a two-dimensional gridded environment containing obstacles. The research has addressed the challenge of determining the most efficient and feasible trajectory. The constraints encompass the initial and final locations of the drone, the locations of randomly placed obstacles, the availability and random placement of battery charging stations within the grid, and the drone's limited travel distance on a single battery charge.

In the model formulation, the drone's single charge travel distance has been set to $2\sqrt{2}$, and one of the constraints is tailored to this assumption. While acknowledging this as a potential limitation, generalizing this parameter is straightforward, and the existing constraints effectively encapsulate the fundamental structure of the generalized constraint.

To tackle this problem, Mixed Integer Programming (MIP) has been employed. The approach primarily centers on the geometric aspects of path planning without introducing time as a variable. The central objective of this research has been to optimize the trajectory, ensuring it is both feasible and of minimal length, thereby facilitating efficient drone navigation in the presence of obstacles. This study offers valuable insights into the geometric considerations of drone path planning, which have direct implications for applications spanning surveillance, logistics, and robotics.

\subsection{Optimal Scheduling for Drone Delivery}
This research has addressed a logistical optimization problem within a two-dimensional grid framework. The primary objective has been to ascertain the optimal trajectories for two drones, each originating from distinct warehouses, to efficiently fulfill multiple delivery requests, each characterized by its unique demand.

This mathematical model has made several key assumptions. Firstly, each drone has been assumed to have a single-package capacity. Warehouses have been presumed to possess an infinite supply of packages. The packages have been assumed to be identical in nature. Additionally, the time required for a drone to move from one grid point to another has been considered to remain constant, irrespective of whether the movement is diagonal or parallel to a side of the square. While each drone commences its journey from a separate warehouse, it possesses the flexibility to revisit any warehouse later for package collection. Importantly, it has been assumed that the drones can navigate the grid without encountering battery limitations. This specific assumption has been consciously adopted to streamline our approach to address the scheduling problem under examination.

To tackle this problem, the research has leveraged Mixed Integer Programming (MIP). The central objective of this model has been to minimize the overall time required to satisfy the delivery demands emanating from various delivery locations, which inherently leads to a reduction in the cumulative travel distance of both drones. To achieve this, a binary variable has been introduced to represent the state of each drone, indicating whether it is engaged in delivery or in transit back to the warehouse. The dynamics of this binary variable are governed by a set of constraints, which in turn guide the satisfaction of supply and demand requirements both at the warehouses and the delivery locations.

This study contributes to the field of logistics by offering a comprehensive approach to optimizing drone-based delivery systems, considering a range of real-world constraints, and employing MIP techniques to minimize delivery time.

\section{Mathematical Formulation}

\subsection{Optimizing the Location of Battery Charging Stations}

\subsubsection{Parameters}
\begin{itemize}

\item O: This is an array storing obstacle locations.

\item $d_{max}$: This is the maximum distance traversable by a drone in a single battery charge.

\item D: This is a distance matrix that gives the shortest distance between any two cells on the grid. It is equal to the Euclidean distance the drone has to cover while taking the shortest path between them considering the fact that obstacles are not traversable.

The distance matrix has been calculated for each pair of neighboring cells in our grid using Dijkstra's algorithm, which is explained in the problem setting above.

For two grid cells (i, j) and (k, l), the distance between them is given by D(i, j, k, l).

\item $B_{depot}$: This traversable (non-obstacle) grid cell, selected at random, serves as the location where a battery charging station is initially placed. This station is essential for implementation in the Miller-Tucker-Zemlin subtour elimination method.

\end{itemize}

\subsubsection{Decision Variables}
\begin{itemize}

\item B: This is a 10 x 10 binary variable such that:
 \begin{equation}
    B(i,\ j) =
    \begin{cases}
      1, & \text{if (i, j) has a battery charging station}\\
      0, & \text{otherwise}
    \end{cases}
  \end{equation}

\item E: This is a 10 x 10 x 10 x 10 binary variable which represents the edges from each point on the grid to every other point on the grid. This variable has been used only for battery charging stations in this problem.

\item T: This is a 10 x 10 variable which has been bounded between 1 and 20 for our case. It represents the time at which a battery charging station is visited. This variable is needed for the Miller-Tucker-Zemlin method for subtour elimination.

\end{itemize}

\subsubsection{Objective} 
The objective is to minimize the total number of battery charging stations.
$$min \sum_{i, j = 1} ^ {10} B(i,\ j)$$

\subsubsection{Constraints}

\begin{itemize}

\item Obstacles: A battery charging station cannot be placed on an obstacle\\
\text{$\forall$ grid cells (i, j) $\in$ O:}
\begin{equation}
    B(i,\ j) = 0   
\end{equation}

\item Reachability: There should be at least one battery charging station within a distance of $d_{max}$ from each traversable grid cell\\
\text{$\forall$ grid cells (i, j) $\notin$ O:}
\begin{equation}
    \sum_{k, l = 1} ^ {10} B(k,\ l) \geq 1,\ \text{$\ni$ D(i, j, k, l) $\leq d_{max}$}    
\end{equation}

\item Edges: A constraint encapsulating the definition of the edges E between grid cells as defined in the decision variables section\\
\text{$\forall$ grid cells (i, j) \& (k, l):}
\begin{equation}
    E(i,\ j,\ k,\ l) =
    \begin{cases}
      1, & \text{if B(i, j) = B(k, l) = 1 \&}\\ & \text{D(i, j, k, l) $\leq$ $d_{max}$}\\
      0, & \text{otherwise}
    \end{cases}
\end{equation}

\item Path Connectedness - Outflow: There should be only one edge emerging from each battery charging station\\
\text{$\forall$ battery charging stations (i, j):}
\begin{equation}
    \sum_{k, l = 1} ^ {10} E(i,\ j,\ k,\ l) = 1   
\end{equation}

\item Path Connectedness - Inflow: There should be only one edge entering each battery charging station\\
\text{$\forall$ battery charging stations (i, j):}
\begin{equation}
    \sum_{k, l = 1} ^ {10} E(k,\ l,\ i,\ j) = 1   
\end{equation}

\item Depot Battery Charging Station: This randomly selected grid cell should contain a battery charging station\\
\text{$\forall$ (i, j) = $B_{depot}$:}
\begin{equation}
    B(i,\ j) = 1
\end{equation}

\item MTZ Subtour Elimination: The visited time variables for each battery charging station, except the randomly selected $B_{depot}$, should be monotonically increasing to eliminate subtours\\
\text{$\forall$ grid cells (i, j) $\neq B_{depot}$ \& (k, l) $\neq B_{depot}$:}
\begin{equation}
    T(k,\ l) \geq T(i,\ j) - 100 * \bigl(1 - E(i,\ j,\ k,\ l)\bigl) + 1
\end{equation}

\end{itemize}

\subsection{Shortest Path Problem}

\subsubsection{Parameters}
\begin{itemize}

\item D: This matrix represents Euclidean distances between grid cells in a 4D space (10 x 10 x 10 x 10). It focuses only on the immediate neighbors of point (i, j) for practicality, as these are the primary points of interest.

\item R: This is an array storing initial \& final points.

\item O: This is an array storing obstacle locations.

\item C: This is an array storing location of charging stations.

\item B: This is a 10 x 10 matrix storing 1's in charging locations:
  \begin{equation}
    B(i, \ j) =
    \begin{cases}
      1, & \text{if (i, j) $\in$ C}  \\
      0, & \text{otherwise}
    \end{cases}
  \end{equation} 

\end{itemize}

\subsubsection{Decision Variables}
\begin{itemize}

\item X: This is a 10 x 10 binary variable such that:
  \begin{equation}
    X(i,\ j) =
    \begin{cases}
      1, & \text{if (i, j) is visited by drone}  \\
      0, & \text{otherwise}
    \end{cases}
  \end{equation}

\item V: This a 10 x 10 matrix storing 1's in visited battery locations:
\begin{equation}
    V = B \odot X \ (\text{$\odot$ := Element-wise multiplication})
\end{equation}

\item E: This a 4D matrix storing active edges:
\begin{equation}
    E(i,\ j,\ k,\ l) = X(i,\ j) * X(k,\ l)
\end{equation}
where (k, l) is a direct-neighbor of (i, j).

\end{itemize}

\subsubsection{Objective} 
The objective is to minimize the total traveling cost of the drone (i.e. total traveling distance).
\begin{equation}
    0.5 * \min \sum_{i, j, k, l=1} ^ {10} D(i,\ j,\ k,\ l) * X(i,\ j) * X(k,\ l)  
\end{equation}

\subsubsection{Constraints}
\begin{itemize}

\item Initial \& Final Points: The initial and final points should be initialized to 1\\
\text{$\forall$ grid cells (i, j):}
\begin{equation}
    X(i,\ j) = 1, \ \text{if (i, j) $\in$ R}
\end{equation}

\item Obstacles: This ensures that a blocked cell cannot be traversed\\
\text{$\forall$ grid cells (i, j):}
\begin{equation}
    X(i,\ j) = 0, \ \text{if (i, j) $\in$ O}    
\end{equation}

\item Visited Battery Locations: A constraint encapsulating the definition of the matrix V as defined in the decision variables section\\
\text{$\forall$ grid cells (i, j):}
\begin{equation}
    V(i,\ j) = B(i,\ j) * X(i,\ j)   
\end{equation}

\item Path Connectedness: Basic flow conservation would suffice to ensure the path should be connected between the initial and final points\\
\text{$\forall$ grid cells (i, j):}
  \begin{equation}
    \sum_{k, l = 1}^{10} E(i,\ j,\ k,\ l) =
    \begin{cases}
      1*X(i,\ j), & \text{if (i, j) $\in$ R} \\
      2*X(i,\ j), & \text{otherwise}
    \end{cases}
  \end{equation}

\item Battery Constraints: For each visited node, there should at least be one battery station accessible for the drone to go to and another one for the drone to come from\\
\text{$\forall$ grid cells (i, j):}
\begin{equation}
    \sum_{k = i-1}^{i+1} \sum_{l = j-1}^{j+1} V(k,\ l) \geq 
    \bigl(1 + X(i,\ j) - V(i,\ j)\bigl) * X(i,\ j)
\end{equation} 
\end{itemize}

\subsection{Optimal Scheduling for Drone Delivery}

\subsubsection{Parameters}

\begin{itemize}
    \item D: This matrix represents Euclidean distances between grid cells in a 4D space (10 x 10 x 10 x 10). It focuses only on the immediate neighbors of point (i, j) for practicality, as these are the primary points of interest.
    \vspace{2pt}
    \item W: This is an array storing warehouse locations.
    \vspace{2pt}
    \item L: This is an array storing delivery locations.
    \vspace{2pt}
    \item Q: This is an array storing delivery demands of respective delivery locations.
\end{itemize}

\subsubsection{Decision Variables}

\begin{itemize}

\item X: This is a 4D binary variable storing the state of the grid:
  \begin{equation}
    X(n,\ i,\ j,\ t) =
    \begin{cases}
        1, & \text{if drone n is present at (i, j)}\\ & \text{at time t}\\
        0, & \text{otherwise}
    \end{cases}
  \end{equation} 

\item f: This is a 2D binary flag variable which specifies the delivery state of the drone:
  \begin{equation}
    f(n,\ t) =
    \begin{cases}
      0, & \text{if drone n is out for delivery}\\ & \text{at time t}\\
      1, & \text{if drone n is coming to warehouse}\\ & \text{at time t}
    \end{cases}
  \end{equation} 

\item T: This is a 1D integer variable storing the total time taken to deliver by each drone.

\end{itemize}

\subsubsection{Objective} 
The objective is to minimize the total traveling cost of the drones (i.e. total traveling distance) and the sum of the individual delivery time of each drone.
\begin{equation}
    \min \Biggl(0.7 * \sum_{n, i, j, k, l, t} d(n, i, j, k, l, t) + 0.3 * \sum_{n} T(n)\Biggl)
\end{equation}
where $d = D(i,\ j,\ k,\ l) * X(n,\ i,\ j,\ t) * X(n,\ k,\ l,\ t+1)$

\vspace{1mm}

\subsubsection{Constraints}

\begin{itemize}

\item Initial Points: The drones should be initialized at their respective warehouse locations\\
\text{$\forall$ n:}
\begin{equation}
    X(n,\ i_{n},\ j_{n},\ 0) = 1,  \ \text{if $(i_{n},\ j_{n})$ $\in$ W}
\end{equation}

\item Path Connectedness: The path formed by the drone should be continuous and connected\\
\text{$\forall$ drones n, grid cells (i, j) \& time t:}
\begin{equation}
    \sum_{k = i-1}^{i+1} \sum_{l = j-1}^{j+1} X(n,\ k,\ l,\ t+1) \geq X(n,\ i,\ j,\ t)
\end{equation}

\item Sequential Movement: At a particular time, each drone should be present at only one node in the grid\\
\text{$\forall$ drones n \& time t:}
\begin{equation}
    \sum_{i, j = 1}^{10} X(n,\ i,\ j,\ t) = 1
\end{equation}

\item Initial Flag: Initialize flag to 0 for all drones\\
\text{$\forall$ drones n:}
\begin{equation}
    f(n,\ 0) = 1
\end{equation}

\item Flag Operation - Delivery: Set flag to 0 when drones are out for delivery\\
\text{$\forall$ drones n \& time t:}
\begin{equation}
    \bigl(1 - f(n, t)\bigl) * \Biggl(f(n, t+1) - \sum_{(i, j) \in L} X(n, i, j, t+1)\Biggl) = 0
\end{equation}

\item Flag Operation - Return: Set flag to 1 when drones are returning to warehouse\\
\text{$\forall$ drones n \& time t:}
\begin{equation}
    f(n, t) * \bigl(f(n, t+1) + X(n, i_{n}, j_{n}, t+1) - 1\bigl) = 0
\end{equation}
\hspace{58mm}if $(i_{n}, j_{n}) \in$ W

\item Supply Satisfaction: Total warehouse supply should precisely match the collective net demand of all delivery locations
\begin{equation}
    \begin{aligned}
        \sum_{(i, j) \in W} \sum_{n, t} X(n,\ i,\ j,\ t) * f(n, t-1)* \bigl(1 - f(n, t)\bigl)\\ = \sum_{(i, j) \in L} Q(i,\ j)
    \end{aligned}
\end{equation}    

\item Demand Satisfaction: The demand of each delivery location should be satisfied\\
\text{$\forall$ (i, j) $\in$ L:}
\begin{equation}
    \sum_{n, t} X(n,\ i,\ j,\ t) * f(n,\ t) * \bigl(1 - f(n,\ t-1)\bigl) = Q(i,\ j)
\end{equation}

\item Collision Avoidance: No grid cell should have more than one drone present simultaneously\\
\text{$\forall$ drone combinations (m, n), grid cells (i, j) \& time t:}
\begin{equation}
    X(m,\ i,\ j,\ t) * X(n,\ i,\ j,\ t) = 0
\end{equation}

\item Delivery Time: This constraint defines the total delivery time for each drone which is used in the objective function\\
\text{$\forall$ drones n:}
\begin{equation}
    T(n) = \max_{t}\ t * f(n,\ t-1) * \bigl(1 - f(n,\ t)\bigl)
\end{equation}

\end{itemize}

\section{Results}

\subsection{Optimizing the Location of Battery Charging Stations}

The only variables in this problem are the location of the obstacles in the grid, and the value of the distance $d_{max}$ which is the maximum distance traversable by the drone in a single battery charge. Note that this distance can be easily correlated with the number of hops the drone can take in a single charge. The initial charging station $B_{depot}$ to be used for the Miller-Tucker-Zemlin subtour elimination, has been strategically positioned within the interior of the grid by random selection, on any of the accessible grid cells devoid of obstacles.

The results have been shown for two different cases. In both cases, $d_{max} = 2\sqrt{2}$ (which corresponds to a maximum of two hops by the drone in any direction), but the number and location of the obstacles have been varied.

For clarity and ease of interpretation, a color-coded grid system has been implemented, using the following key:
\begin{itemize}
    \item Black - Location of obstacles not traversable by the drone (predetermined for a problem setting)
    \item Blue - Location of battery charging stations (given by the model)
    \item White - Grid points traversable by the drone
\end{itemize}

Note that the battery charging stations given by the model have been placed only on the traversable grid points, which were previously white.

\begin{figure}[H]
\centerline{\includegraphics[width=0.29\textwidth]{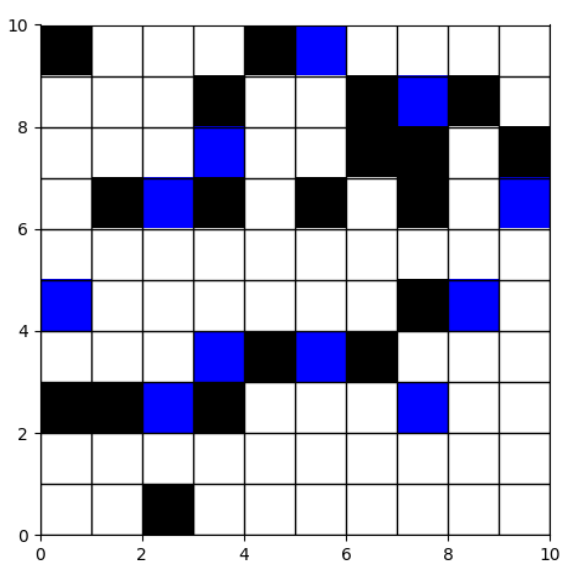}}
\caption{Minimum number of battery charging stations = 11}
\label{fig:1}
\end{figure}

\begin{figure}[H]
\centerline{\includegraphics[width=0.29\textwidth]{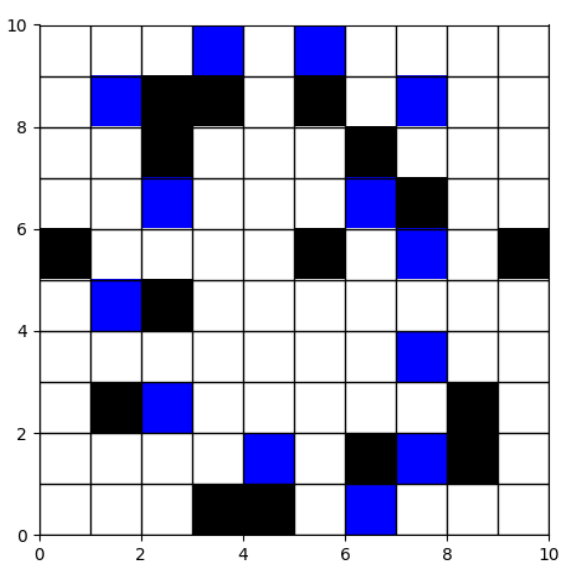}}
\caption{Minimum number of battery charging stations = 13}
\label{fig:2}
\end{figure}

Figures \ref{fig:1} and \ref{fig:2} show that the minimum number of battery charging stations required to service the entire grid changes with a change in the number and location of the obstacles. The selection of the initial $B_{depot}$ significantly influences the result. An apt placement can potentially reduce the required number of battery charging stations.

It can also be verified that each and every grid cell that does not contain an obstacle has at least one battery charging station within a distance of $2\sqrt{2}$. This ensures that each and every node is reachable from a battery charging station. Secondly, all the battery charging stations are themselves connected to each other, i.e. there exists a path from each battery charging station to every other battery charging station. This directly implies that there exists a path from each free grid cell to every other free grid cell.

\subsection{Shortest Path Problem}

In this problem, the locations of the battery charging stations, obstacles, and the initial and final points are randomly chosen and act as inputs to the model. The initial and final points can be considered to be the warehouse and delivery locations in the real-world scenario. The model described above gives the shortest feasible connected path, keeping in mind the battery constraints. If the obstacles and battery charging stations have been distributed in such a way that a path cannot exist, then the computation is terminated prematurely and the model just returns an infeasible solution.

For clarity and ease of interpretation, a color-coded grid system has been implemented, using the following key:
\begin{itemize}
    \item Black - Location of obstacles not traversable by the drone (predetermined for a problem setting)
    \item Blue - Location of battery charging stations (predetermined for a problem setting)
    \item Red - Initial and final points of the drone (predetermined for a problem setting)
    \item Purple - Location of battery charging stations that are visited by the drone (given by the model)
    \item Green - Grid points that are visited by the drone (given by the model)
    \item White - Other grid points
\end{itemize}

Two cases of this problem with different inputs to the model have been shown, which give different shortest paths to be followed by the drone in accordance with the battery limitations.

\begin{figure}[H]
\centerline{\includegraphics[width=0.30\textwidth]{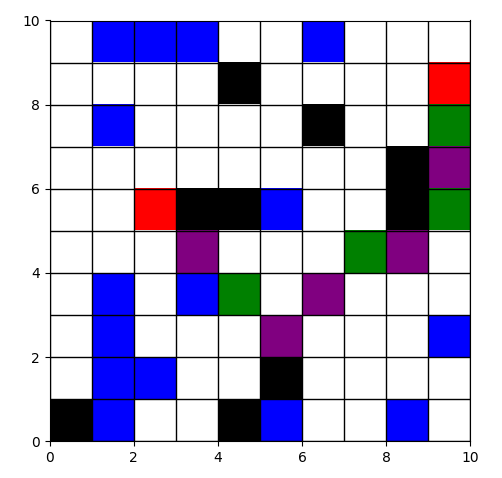}}
\caption{Shortest path for drone (1)}
\label{fig:3}
\end{figure}

\begin{figure}[H]
\centerline{\includegraphics[width=0.30\textwidth]{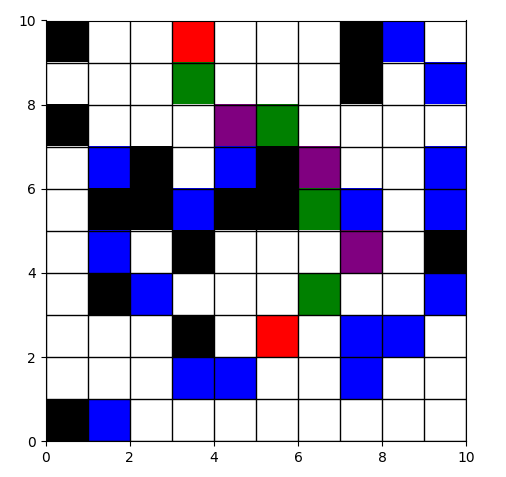}}
\caption{Shortest path for drone (2)}
\label{fig:4}
\end{figure}

In figures \ref{fig:3} and \ref{fig:4}, it can be seen that for the given initial and final points of the drone, the shortest path is pretty straightforward if not for the battery constraints. The battery limitation forces the drone to move in a convoluted manner, thus ensuring that its battery does not die out during its movement. It can be verified that these are the shortest feasible paths for these problem settings.

\begin{figure*}[t]
\centerline{\includegraphics[width=1.1\textwidth]{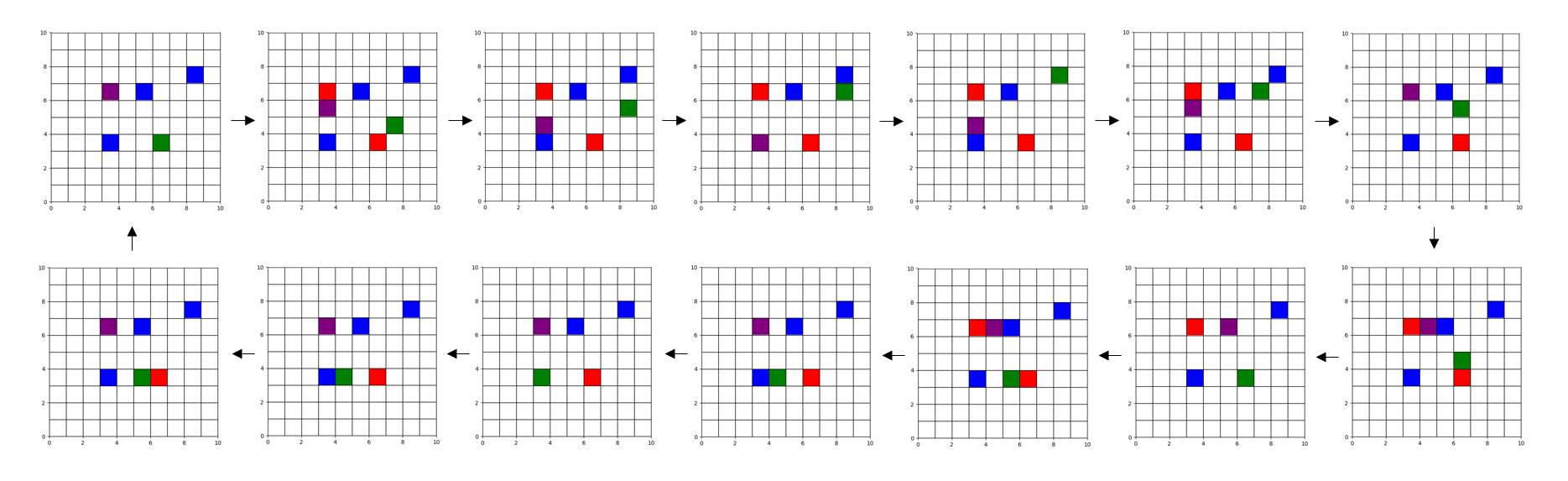}}
\caption{Optimal Scheduling}
\label{fig:5}
\end{figure*}

\subsection{Optimal Scheduling for Drone Delivery}

This problem takes as input the warehouse locations, the delivery locations, the number of drones, and the demand of each delivery location, and outputs the position of the drones at each time instant.

For clarity and ease of interpretation, a color-coded grid system has been implemented, using the following key:
\begin{itemize}
    \item Red - Warehouse locations for the drones (predetermined for a problem setting)
    \item Blue - Delivery locations (predetermined for a problem setting)
    \item Purple - Location of drone 1 (given by the model)
    \item Green - Location of drone 2 (given by the model)
    \item White - Other grid points
\end{itemize}

Figure \ref{fig:5} depicts the color coding and the movement of the drones in a total of fourteen timesteps. In this case, two warehouse locations for the drones and three delivery locations have been provided. The delivery locations have varied demands for packages, which are tended to by the two drones. The delivery location on the bottom left has a demand of two, while the other two delivery locations have a demand of one. The drones start their deliveries at the warehouses, go to the delivery locations to deliver the packages and come back to the warehouses at the end. Furthermore, to implement collision avoidance, two drones cannot be at the same grid point at the same time, be it a warehouse, delivery location, or free space. But the drones are free to choose at the end of each delivery, which warehouse they would like to visit to pick up the package from for the next delivery.

Since the total time taken by the drones to deliver all the packages has been minimized, we see that the drones choose which delivery locations they would visit and start their journey simultaneously. The delivery location to the top right is attended by the green drone, the delivery location in the middle by the purple drone, and the bottom left delivery location is visited by both the purple and green drones since it has a delivery demand of two.

\section{Conclusions}

This study tackles three significant challenges associated with the implementation of drone delivery systems: (1) the optimal positioning of battery charging stations, (2) solving the shortest path problem for drones between their initial and final locations constraint to their single battery charge travel distance, and (3) efficiently scheduling multiple drones amidst numerous warehouses and delivery locations, each with varied demands.

This study represents the environment as a 2D grid incorporating obstacles. This simplification was chosen for its practicality, as top-view images of maps can be pixelated and treated as grids. This methodology provides a solid foundation, offering the potential for extension into a 3D grid, where any structure can be modeled as an obstacle in the three-dimensional space and a battery charging station can be placed either on the ground, terrace or a selection of floors like the refuge areas in a building.

Our approach to placing battery charging stations involves applying the Miller-Tucker-Zemlin subtour elimination method to prevent the formation of clusters. This method necessitates the initialization of a depot node. In a grid with specified obstacle configurations, the ultimate positioning of battery charging stations is notably influenced by the chosen location of this depot node, even though the overall count of charging stations remains relatively constant.

For the case of optimal scheduling, the chosen objective function incorporates both a distance metric and the cumulative time required for the completion of all deliveries by the drones. While distance and time are typically directly proportional, merely minimizing one variable wouldn't suffice in our context. This arises due to the consideration of diagonal and lateral hops as temporally equivalent in our model, despite their differing distances. Considering only time in the objective function results in unnecessary zig-zag motion due to the temporal equivalence of both diagonal and lateral movements. On the other hand, minimizing only distance can lead to drones hovering in certain solutions due to the absence of restrictions on delivery time.  

As this is one of the first papers comprehensively addressing drone delivery systems, many potential areas for future research arise. One avenue to explore initially involves combining the three cases and formulating a holistic solution. In essence, when presented with a geographical area, potential warehouse locations, and delivery points, the model would strategically position battery charging stations before tackling the scheduling problem for multiple drone deliveries, encompassing the shortest path challenge. An intriguing extension involves venturing into three-dimensional space, paving the way for integrated ground and air transport in delivery systems. Another challenging aspect worth exploring pertains to bi-level optimization, where the optimal placement of battery charging stations and the determination of the shortest path become interdependent. The resolution of the former has the potential to adversely influence the solution to the latter, adding complexity to the optimization process.

\section{Acknowledgements}

The authors would like to thank Prof. Avinash Bhardwaj from the department of Industrial Engineering and Operations Research at the Indian Institute of Technology Bombay for his technical support and guidance throughout the research process.

\bibliographystyle{plain}
\nocite{*}
\bibliography{drone}

\end{document}